\documentclass[12pt,a4paper]{article}
\usepackage{amsmath, amssymb, theorem, latexsym}
\usepackage{graphicx}
\usepackage{color}
\newtheorem{theorem}{Theorem}[section]
\newtheorem{lemma}[theorem]{Lemma}
\newtheorem{proposition}[theorem]{Proposition}

\newtheorem{remark}[theorem]{Remark}

\allowdisplaybreaks

\newcommand\finbox{~\hfill$\Box$}%

\def\la {{\lambda}}

\def \Inte{{\rm Int\,}}

\newcommand {\pa}{\partial}

\numberwithin{equation}{section}
\title{Minimal partitions for anisotropic tori}

\author{ B. Helffer (Universit\'e Paris-Sud 11) \\
 and\\
 and T. Hoffmann-Ostenhof (University of Vienna)}
 \date{}

\begin{document}
\maketitle

\begin{abstract}
We analyze spectral minimal $k$-partitions for the torus. In continuation with what we have obtained for thin annuli  or thin strips on a cylinder (Neumann case), we get similar results for anisotropic tori.\\


\end{abstract}

\section{Introduction}\label{section1}
For $a\geq b>0$,  consider the Laplacian on the $2D$-torus : $T(a,b):= \mathbb S^1 (\frac{a}{2\pi}) \times \mathbb S^1 (\frac{b}{2\pi}) \,$.  Concretely, we can also consider
\begin{equation}\label{defR}
\mathcal R(a,b)=(0,a)\times (0,b),\,
\end{equation}
and the Laplacian on $\mathcal R(a,b)$ with periodic boundary conditions but except for the pictures this is not the most convenient point of view. It is indeed better to think of the torus as a compact regular manifold.\\
 We can, following \cite{HHOT},   consider $k$-partitions $\mathcal D$ of the torus, i.e. families of disjoint open sets $(D_1,\dots,D_k)$ of the torus
  and the  sequence  of partition energies $\mathfrak L_k(T(a,b))$ obtained by minimizing over $\mathcal D$  of the torus some  energy  defined by 
   \begin{equation}\label{energ}
 \Lambda_k (\mathcal D) = \max_j \lambda (D_j)\,,
 \end{equation}
 where $\lambda(D_j)$ is the ground state energy of the Dirichlet Laplacian in $D_j$.  We then define
 \begin{equation}
 \mathfrak L_k(T(a,b) ) := \inf_\mathcal D \Lambda_k(\mathcal D)\,,
 \end{equation}
 where the infimum is over all the $k$-partitions of $T(a,b)$.
  A minimal $k$-partition is a partition whose energy is $ \mathfrak L_k(T(a,b) ) $\,.
 As in the case of an open set in $\mathbb R^2$, minimal $k$-partitions exist and are strong and  regular (see Section \ref{section2}).  Without loss of generality, we consider the case $a=1$. Note that for the torus, when $b<1$, $\la_1=0$ and that $\la_2=\la_3=4 \pi^2 $. Hence $\mathfrak L_3>\la_3$ and using the results of \cite{HHOT} (extended to the case of the torus) (see Theorem \ref{thHHOT} in Section \ref{section2}) the associated 
minimal $3$-partition cannot be nodal, i.e a partition obtained as the nodal domains of an eigenfunction. 
 On the other hand for $k=4$, we see that $\lambda_4 = 16 \pi^2$ for $b <\frac 12 $ and that any corresponding eigenfunction has four nodal domains. So the minimal $4$-partition is nodal.  Our aim in this paper is to describe what are the minimal $k$-partitions. Our main result is the following:
\begin{theorem}\label{torus}~\\
There exists $b_k >0$ such that, if 
 $b<b_k$,  $\mathfrak L_k(T(1,b))= k^2 \pi^2$ and the corresponding minimal $k$-partition 
$\mathcal D_k=(D_1,\dots ,D_k)$ is  represented in $\overline{\mathcal R(1,b)}$   by 
\begin{equation}\label{minD}
D_i= (\, (i-1)/k\,,\,i/k\,)\times [0\,,\, b\,)\,,\mbox{ for } i=1,\dots, k\,.
\end{equation}
Moreover we can take  $b_k =  \frac{2}{ k}$ for $k$ even and $b_k=\frac 1k$ for $k$ odd.
\end{theorem}

\begin{figure}[h!]\label{Fig1}
\begin{center}
   \caption{One candidate for the minimal $3$-partition represented in $\mathcal R(1,b)$.}\label{figa} 
 \includegraphics[height=7cm]{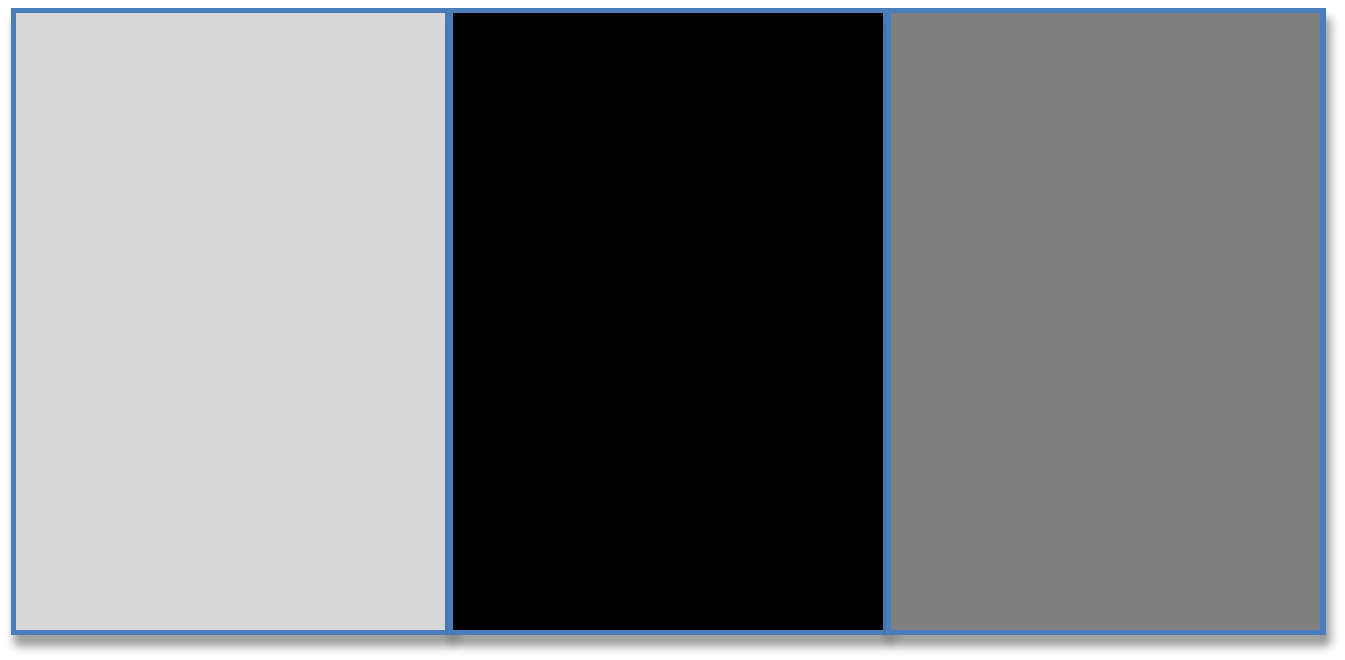}
\end{center}
\end{figure}

Note  that the boundaries of the $D_i$ in $T(1,b)$ are just $k$ circles (see Figure~\ref{Fig1} where these circles are represented by vertical segments). 

\begin{remark}~\\
This result is a complement to  what we have obtained for thin annuli or strips on a cylinder (in the case of the Neumann condition) \cite{HHO2}. Its proof requires new ideas which hopefully can be used for other compact surfaces. We recall that the case of thin annuli with  Dirichlet conditions  is still open ($k$ odd). 
For minimal $k$-partitions of the torus, we will  at the end of Section \ref{section2} prove  that the statement of the theorem holds for $k$ even (the minimal partitions are nodal)  and $b_k=\frac 2 k$ cannot be improved (see also Section \ref{Section7} for further discussion). For $k$ even and  $\frac 2k<b < \frac{2}{k-2}$, the $k$-th eigenfunction does not have $k$ nodal domains. Hence  it remains to give the proof  of our theorem  for $k$ odd ($k\geq 3$).\\
We also recall that in the case $k=3$, the problem was solved in \cite{HHOT:2010} for the sphere $S^2$ and is still open for the disk \cite{HHOT:2010} and the square \cite{BHHO}.
\end{remark}

\section{Reminder on the properties of minimal partitions}\label{section2}
Let us first recall in more detail the properties of minimal $k$-partitions. 
The notion of minimal partition  was first introduced for an open set $\Omega$ in $\mathbb R^2$ in \cite{HHOT} (see references therein). We just present the corresponding definitions for the torus (or more generally on a compact Riemannian manifold).
We recall that a $k$-partition on the torus is simply a family $\mathcal D$ of $k$-disjoint open sets $(D_i)_{i=1,\dots,k}$. Such a  partition is called {\bf strong} if 
$\cup \overline{D_i} = T(1,b)$ and $\Inte ( \overline {D_i}) = D_i$ for any $i$.
Attached to a strong  partition, we  associate a closed
set in $ T(1,b)$, which is called the {\bf boundary set}  of the partition~:
\begin{equation}\label{assclset} 
N(\mathcal D) = \overline{ \cup_i  \partial D_i }\;.
\end{equation}
$ N(\mathcal D)$ plays the role
 of the nodal set (in the case of a nodal partition).  We have recalled in the introduction the notion of minimal $k$-partitions. 
 As in the case of an open set in $\mathbb R^2$, minimal $k$-partitions exist and are strong and  {\bf regular}  in the following sense.
We call a partition  $\mathcal D$ regular if its associated
 boundary set  $ N(\mathcal D) $, has the following properties~:\\
(i)
Except for finitely many distinct $ x_i\in  N$
 in the neighborhood of which $ N$ is the union of $\nu_i= \nu(x_i)$
smooth curves ($ \nu_i\geq 3$) with one end at $ x_i$,  
$ N$ is locally diffeomorphic to a regular 
curve.\\
(ii) $ N$  has the {\bf equal angle
  meeting
 property}.
The $x_i$ are called the critical points and define the set
$X(N)$. 
By {\bf equal angle meeting property}, we mean that   the half curves meet with equal angle at each critical
 point of $ N$.\\
  In the case of an open set we have also points $y_j$ at the boundary and we call this set $Y(N)$.\\
  
 Moreover, the minimal $k$-partitions are bipartite, i.e.  can be colored by two colors (neighboring domains  have different colors),  if and only if they are nodal (i.e. corresponding to the nodal domains of an eigenfunction of the Laplace-Beltrami operator). Another important statement established in \cite{HHOT}  is:
  \begin{theorem}\label{thHHOT}~\\
  A $k$-partition  consisting  of   the $k$ nodal domains of an eigenfunction corresponding to the $k$-th eigenvalue $\lambda_k$ of the Laplacian is a minimal $k$-partition.
   \end{theorem}
In general  one could just say that by the well known  Courant  nodal theorem the  number of nodal domains of an eigenfunction $u_k$  associated with $\lambda_k$ is at most $k$.   The eigenpair  $(u_k,\lambda_k)$ is called {\bf Courant sharp} if the number of nodal domains is exactly $k$. Theorem \ref{thHHOT} is moreover optimal 
as has  been  proven in \cite{HHOT}:
\begin{theorem}\label{thHHOT2}~\\
 A  nodal minimal  $k$-partition  corresponds necessarily to a Courant sharp pair.\end{theorem}

\paragraph{First application: proof of Theorem \ref{torus} in the even case.}~\\
For the torus $T(c,d)$, the eigenvalues are given by $4 \pi^2 (\frac{m^2}{c^2} + \frac{n^2}{d^2}) $ ($(m,n)\in \mathbb N^2$ where $\mathbb N$ denotes the set of the non-negative integers) with a corresponding basis 
 given by 
 \begin{itemize}
 \item
 $(x,y)\mapsto \cos (2\pi m \frac xc)\cos (2\pi n \frac y d)$, 
 \item
 $(x,y)\mapsto \cos (2\pi m \frac xc)\sin (2\pi n \frac y d)$, 
 \item
 $(x,y)\mapsto \sin  (2\pi m \frac xc)\cos (2\pi n \frac y d)$
 \item
  and $(x,y) \mapsto \sin(2\pi m \frac xc)\sin  (2\pi n \frac y d)$ 
  \end{itemize}
  (with suitable changes when $m$ or $n$ vanishes). For example, for $n=0$, we get $(x,y)\mapsto 1$ for $m=0$ and
   $(x,y)\mapsto  \cos (2\pi m \frac xc)$ and  $(x,y)\mapsto  \sin (2\pi m \frac xc)$ for $m>0$. These eigenfunctions have $(2m)$ nodal domains on the torus.  When $k$ is even and $k <  \frac{2c} {d} $, we get the existence of an $k$-th eigenfunction with exactly $k$ nodal domains (corresponding to $m=\frac k 2$ and $n=0$).   What will be important in our
problem is that
Theorem  \ref{thHHOT}  implies that for $k$ even and $c>d>0$ the minimal
$k$-partitions of $T(c,d))$ are nodal for the case that
$k<2c/d$. The corresponding energy is $\frac{\pi^2k^2}{c^2}$.   Hence we have completed the proof of Theorem \ref{torus} for the even case.\\
 
  We also observe that for $k$ odd ($k>1$) the minimal $k$-partitions cannot be nodal.\\ We will prove  that, when $\frac cd$ is small enough,   the minimal $k$-partitions can be lifted into a Courant sharp $(2k)$-partition on the covering $T(2c,2d)$.The $k$-partition appearing in Theorem \ref{torus} corresponds actually to a nodal partition on this covering and this implies the result.  The existence of this lifting will be proved in Section \ref{section5}.

\section{Necessary conditions}\label{section3}
The computation of  the energy of the   $k$-partition \eqref{minD} leads immediately to the following upper bound for $\mathfrak L_k$.
\begin{proposition}~\\
\begin{equation}\label{en1}
k^2 \pi^2 \min (1,b^{-2}) \geq \mathfrak L_k(T(1,b))\,.
\end{equation}
\end{proposition}

Using this upper-bound we can give  necessary  conditions on $k$-partitions to be minimal.

\begin{proposition}\label{prop3.2}~\\
If  $b < \frac1 k\,,$ there is no minimal $k $-partition $\mathcal D= (D_1,\dots ,D_k)$ of the torus with one $D_i$  homeomorphic to a disk.
\end{proposition}

The proof is by contradiction. 
Let $\mathcal D =(D_1,\dots ,D_k)$  be a 
minimal $k$-partition such that, say $D_1$ is homeomorphic to a disk. 
Then,   the  pullback  $\widehat D_1$ of $D_1$ in the universal covering $\mathbb R^2$ is a union of bounded components $\widehat D_1^{k,\ell}$ (with $(k,\ell)\in \mathbb Z^2$) such that $\widehat D _1^{k,\ell} + (m,nb) = \widehat D_1^{k+m,\ell+n}$.
Moreover $\widehat D_{1}^{0,0}$ has same area as $D_1$ and $\lambda (D_1) = \lambda(\widehat D_1^{0,0})\,$. \\
Looking at a lower bound for  $\lambda(\widehat D_1^{0,0})$, one could first think of using Faber-Krahn's inequality but it is better  to come back to the first step of one  proof of the Faber-Krahn  inequality which is based on the Steiner symmetrization (see for example the book \cite{He} (Section 2.2) or the expository talk  \cite{Tr}).\\
 We now observe that each vertical  slice has a total length less than $b$. We now apply the Steiner symmetrization with respect to the horizontal line $y= \frac b2$. It is immediate to see that the image $S(\widehat D_{1}^{0,0})$ of  $\widehat D_{1}^{0,0}$ is contained in a rectangle $\widehat R_b$ in the form $(-\ell_b,\ell_b)\times (0,b)$ for some $\ell_b >0$. Now it is well known that in this symmetrization we have:
 $$
  \lambda(\widehat D_1^{0,0})\geq \lambda (S(\widehat D_{1}^{0,0}))\,,$$
  and by monotonicity
  $$
  \lambda (S(\widehat D_{1}^{0,0}) ) \geq \lambda (\widehat R_b) =\pi^2 (b^{-2} +\ell_b^{-2})> \pi^2 b^{-2}\,.
  $$
  This leads to
  \begin{equation}\label{St1}
  \lambda(D_1) > \pi^2 b^{-2}\,,
  \end{equation}
  hence, using \eqref{en1}, to
  \begin{equation}\label{St2}
  b > \frac  1k\,.
  \end{equation}
  This gives the contradiction.
 
 \section{Around Euler's formula}
\subsection{Standard Euler's formula}
In the case of an open set $\Omega$ in $\mathbb R^2$, observing that the Euler characteristic  of $\Omega$ is $2$,  we have for a regular minimal $k$-partition $\mathcal D$:
 \begin{equation}\label{Euler0}
 k = {\bf b_1} - {\bf b_0} +1 + \sum_i \left(\frac{\nu(x_i)}{2} -1\right) + \frac 12  \sum_{j}  \rho(y_j)\,.
 \end{equation}
 where ${\bf b_0}$ is the number of components of $\pa \Omega$,  ${\bf b_1}$ is the number of components of $\pa \Omega \cup N$,  $\nu(x_i)$ and $\rho(y_j)$
the numbers of arcs associated with the singular points $x_i\in X(N)$ of the boundary set $N=N(\mathcal D)$   in $\Omega$, respectively  with the points $y_j$  of the boundary set contained in $\pa \Omega$. We denote by $X(N)$ the set of the $x_i$'s and by  $Y(N)$ the set of the $y_j$'s.
\subsection{Euler's formula on the torus and applications} 
In the case of a flat  compact surface $M$ without boundary, it is easier to formulate Euler's formula by using the Euler's characteristics of $M$ and of the elements of the partition
 $\mathcal D =(D_1,\dots,D_k)$. The formula reads 
\begin{equation}\label{Euler1} 
  \sum_\ell \chi( D_\ell) =  \chi (M)  +  \sum_{i}\left(\frac{\nu(x_i)}{2}-1\right) \,,
\end{equation}
and is a direct consequence\footnote{Thanks to P. B\'erard for giving us the reference.} of the Gauss-Bonnet formula applied in each $D_i$ (see for example \cite{SST}).\\
We recall that for the torus: $\chi(T(a,b)) =0$, for the disk $B$:  $\chi (B)=1$, for the annulus $A$:  $\chi (A)=0$ and for the sphere $\mathbb S^2$: $\chi (\mathbb S^2)=2$. Hence, in the case of the torus,  \eqref{Euler1} becomes:
\begin{equation}\label{Euler2}
\sum_{\ell =1}^k \chi(D_\ell) =  \sum_{i}\left(\frac{\nu(x_i)}{2}-1\right) \,.
\end{equation}

\begin{proposition}\label{cor2}~\\
A minimal partition $\mathcal D=(D_1,...,D_k)$ for which no $D_\ell$ is 
homeomorphic to the disk  satisfies
$X(N(\mathcal D))=\emptyset\,$.
 \end{proposition}
 {\bf Proof.}\\
 The assumption implies that $\chi(D_\ell)\leq 0$, for $\ell=1,\dots,k$. Then we immediately get  from \eqref{Euler2} that $\chi (D_\ell) =0$ and that $X(N)=\emptyset$.  

 \section{Lifting argument}\label{section5}

\begin{proposition}\label{lastprop}~\\
Suppose 
$\mathcal D=(D_1,\dots,D_k)$ is a minimal k-partition on the torus $T(1,b)$
for which all the $D_i$ are not homeomorphic to the disk and 
$X(N(\mathcal D))=\emptyset$. Then $\mathcal D$ can be lifted to a 
bipartite
$(2k)$- partition of $T(2,2b)$.
\end{proposition}

The initial guess was that   a double covering will suffice 
but this is not always the case. One can  construct (see Figure \ref{fig1}) a $3$-partition
of the torus without critical point,  for which it is necessary to construct a covering of order $4$,  $T(2,2b)$ of the torus
  (doubling in each direction) in order to get a bipartite $6$-partition (see Figure \ref{fig2}).
  
  {\bf Proof of Proposition \ref{lastprop}.}\\
 One can classify all the possible topological types of these partitions. 
 The $k$ open sets of the partition have the same topological type. Each open set can be deformed by a retraction  onto a simple closed   line without self-intersection.
Hence  the classification  corresponds to the classification of closed  lines on the torus without self-intersection that are not homotopic  to a point (the so-called torus knots). They  correspond (see \cite{Ha}, p.~47, Example 1.24)  to 
 lines generically denoted by $\ell_{p,q}$ turning $p$ times  around one  horizontal circle and $q$ times around the  vertical  one, with $p$ and $q$ mutually prime (except if $q=0$ , $p=1$ or $p=1$, $q=0$).
  Figure \ref{fig1} corresponds to $p=1$, $q=1$.  The candidate for the minimal $3$-partition when $b$ is small corresponds to $p=1$, $q=0$.  Another example is given in the first subfigure of Figure \ref{figxa}, which represents a closed line on the torus with $p=3$ and $q=2$.
  We go to a suitable double covering so that either $p$ or $q$  is multiplied
by $2$; so the greatest common divisor
$D(p,q)=2$. There are two  cases  : $pq$ odd or $pq$ even (with $p$ or $q$ odd). In the first case we choose $T(2,2b)$ and in the second case the minimal choice is $T(1,2b)$ or $T(2,b)$ but $T(2,2b)$ is also suitable, the important point being that $D(2p,2q)=2$. On the covering $T(2,2b)$, the pull-back of our closed line $\ell_{p,q}$ in $T(1,b)$ is the union of  two distinct  closed lines in $T(2,2b)$. Coming back to the $k$-partition, the lifting to $T(2,2b)$ leads to a $(2k)$-partition.
This ends the proof of the proposition.
  \begin{remark}\label{rem5.2}~\\
  When $p$ and $q$ are not mutually prime, our constructions lead, as explained in \cite{Ha}  to $D(p,q)$ connected closed lines, where $D(p,q)$  is the greatest common divisor of $p$ and $q$. The second subfigure  of Figure \ref{figxa}  corresponds to the case  $p=4$ and $q=2$.\\
 To understand the point, take the closure of $\mathcal R(p,q)$  (see \eqref{defR}) and consider the intersection of the  lines of equation
  $y=-x +c $ ($c\in \mathbb Z$) with $\overline{\mathcal R(p,q)}$. If we project on the corresponding torus and look at the number of connected components obtained on the torus, then we observe that this number is 
   $D(p,q)$ (see  the second subfigure  of Figure \ref{figxa}  which has two components). When $D(p,q)=1$, we get a single closed line of the torus. After a suitable dilation, we can then come back to $T(1,b)$.\\
   When $D(p,q)\neq 1$, it is  not possible to find a continuous  closed  line on the torus without  self-intersection with winding pair $(p,q)$.
   \end{remark}
  
  \section{End of the proof of Theorem \ref{torus}}
    We deduce from Propositions \ref{prop3.2}, \ref{cor2} and \ref{lastprop} that, if $b<\frac 1 k$ ($k$ odd), then any minimal $k$-partition can be lifted into a $(2k)$-partition of $T(2,2b)$ with the same energy $\mathfrak L_k(T(1,b))$.
  We need to look at the spectrum of the Laplacian  on  the $4$-covering $T(2,2b)$ and to determine under which condition the $(2k)$-th eigenvalue is Courant sharp.
  The eigenvalues are given by
  $  \pi^2 ( \ell^2 +  m^2/ b^2)$.  If $b<\frac 1k$, the $(2k)-th $ eigenvalue corresponds to $m=0$ and $\ell =k$, and we are in a Courant sharp situation. Theorem \ref{thHHOT} implies that  
  $$
  \pi^2 k^2 = \mathfrak L_k(T(2,2b))\leq \mathfrak L_k(T(1,b)) \,.
  $$
  Having in mind \eqref{en1}, this ends the proof of the theorem in the odd case.
     
   \begin{remark}~\\
The ideas in the proof might lead to results concerning minimal partitions for other
 "thin" compact surfaces.
\end{remark}

    \section{More on the Courant sharpness of eigenfunctions for the case that $b^2$ is irrational.}\label{Section7}
We recall (see after Theorem \ref{thHHOT2}) that on $T(1,b)$, the associated eigenvalues are given by 
\begin{equation}\label{spectruma}
\la_{m,n}(1,b)=4\pi^2(m^2+\frac{n^2}{b^2})\,. 
\end{equation}
If $m,n>0$ and if  $b^2$ is irrational, then we have multiplicity $4$. 
Following some ideas which we presented already in \cite{HHOT} for rectangles and the disk we have the following 
result. 
\begin{theorem}\label{theorem7.1}~\\
 Suppose $b^2$ is irrational. If $\min (m,n) \geq 1$, then there is no Courant sharp pair $(u,\lambda_{m,n})$.
 \end{theorem}
The proof is based on the following 
\begin{proposition}\label{prop7.2}~\\
For $m, n>0$ any eigenfunction $u$ corresponding to $\lambda_{m,n}$ has at most  $4mn$ nodal domains.  Moreover the only eigenfunctions with exactly $4 mn$ nodal domains
 have the form $\cos (2\pi mx +\theta_1) \cos (2\pi n\frac{y}{b} + \theta_2)$ for some constants $\theta_1$ and $\theta_2$. The other eigenfunctions have  $2 D (m,n)$ nodal domains, where $D(m,n)$ is the greatest common divisor of $m$ and $n$.
\end{proposition}
{\bf Proof of the proposition}\\
We first observe that a general eigenfunction associated with $\lambda_{m,n}$ can be written in the form:
\begin{equation}\label{thefu}
u = \mu \left(\cos2\pi  mx \cos (2\pi n\frac{y}{b}+\theta_1) + \lambda \sin 2\pi mx  \cos (2\pi n\frac{y}{b}+\theta_2)\right )\,,
\end{equation} 
with $\mu\neq 0$.\\
Note that it is only here that we use the fact that $b^2$ is irrational. By rotation, we can reduce to the case when $\theta_2=0$ and we write $\theta=\theta_1$.\\
Then after dilation and rotation, the proof is based on the following lemma:
\begin{lemma}~\\
Except when $\lambda=0$ or $\theta\equiv \frac \pi 2 \; mod (\pi)$, the nodal set of the  function $u_{\lambda,\theta} :=\cos 2\pi x \cos (2\pi y +\theta) + \lambda \sin2\pi  x \sin 2\pi  y$ has 
 no critical zero.
\end{lemma}
Let us look at the critical zeroes of this functions. They should satisfy:
\begin{equation}
\begin{array}{ll}
\cos 2\pi  x \cos (2\pi y+\theta) + \lambda \sin2\pi  x \sin2\pi y &= 0\,,\\
-\sin 2\pi x \cos (2\pi y+\theta) + \lambda \cos 2\pi x \sin 2\pi y&=0\,,\\
 -\cos 2\pi x \sin (2\pi y+\theta) + \lambda \sin 2\pi x \cos 2\pi y&=0\,.
\end{array}
\end{equation}
We assume $\lambda \neq 0$. Suppose that this system has a solution. 
The two first equations imply
 $\cos (2\pi y+\theta)=0$ and $\sin 2\pi y =0$. This implies $\cos \theta =0$.\\
 Hence, when $\cos \theta \neq 0$, our function $u_{\lambda,\theta}$ has no critical zero. 
 
 \begin{lemma}~\\
 For $\lambda\neq 0$, $\theta_2=0$, and $\cos \theta \neq 0$, the nodal partition of the function $u$ of \eqref{thefu} has $2 D(m,n)$ components.
 \end{lemma}
 In each connected component of the set $\mathcal A: =\{(\lambda,\theta)\,|\, \lambda \neq 0, \cos \theta \neq 0\}$ in $\mathbb R^2$
  the number of nodal domains is constant. Hence it is enough to determine this number for one specific pair $(\lambda,\theta)$ in each component of $\mathcal A$.
  It is enough to consider $\lambda =\pm 1$ and $\theta \equiv  0\,  ({\rm mod\,} \pi)$, where the computation of the number of nodal domains is immediate (see Remark \ref{rem5.2}) and equal to $2 D(m,n)$.\\
  
  Note that when $\cos \theta =0$, we get a product $$ u_{\lambda,\theta}:= \sin 2\pi n y (\lambda \sin 2\pi m x \pm \cos 2\pi m x)$$ which has $4mn$  nodal domains.

\begin{remark}~\\
This is not clear for the case that $b^2$ is rational, since then higher multiplicities 
could  occur and we do not know how to exclude the possibility of a higher number of nodal domains 
in higher dimensional eigenspaces. 
\end{remark}
{\bf Proof of Theorem \ref{theorem7.1}}\\
We give two alternative proofs (the second is geometric and inspired by arguments developed in \cite{HHOT}):
\paragraph{Proof 1}
If $\inf(n,m)\geq 1$, then $\lambda_{m,n} = \lambda_{k(n,m)}$
 with  $k(m,n)\geq 4 mn + 2m + 2n -2$. This is obtained by just adding the
multiplicities of the eigenvalues $\lambda_{m',n'}$
with $m'\leq m$, $n'\leq n$, $(m',n')\neq (m,n)$. On the other hand, Proposition \ref{prop7.2} says that any eigenfunction has at most $4 mn$ domains (if $\inf (m,n)\geq 1$). Hence it cannot be Courant sharp.\\

\paragraph{Proof 2}
According to Proposition \ref{prop7.2}  it is enough to consider eigenfunctions in the form $\sin (2\pi mx +\theta_1) \sin (2\pi n\frac{y}{b} + \theta_2)$ and to show that it cannot correspond to a Courant sharp case.
Consider for simplicity  the situation that $m=1=n$ and $\theta_2=\theta_1=0$. Then (up to a rotation) the eigenfunction is given by $u_{1,1}=\sin 2\pi x \sin( 2\pi y/b)$. 
The zeros are given by the zeros of the sines. In particular we can for instance consider the zero 
given by $y=b/2$ and $y =0$. Consider the $P_1=\{(x,y)\in T(1,b)\:| 0<y<b/2\}$ and $P_2=\{(x,y)\in T(1,b)\:b/2<y<b\}$
and consider $N_i=\{(x,y)\in \overline P_i\:|(x,y)\in N(u_{1,1}(x,y))\}\,$, where $N(u)$ denotes the zeroset of $u$. 
Suppose we have a minimal partition corresponding to this eigenfunction. Then we can rotate for instance 
$N_1\,$,  so that the zeros $x=0\,,\,x=1/2$ are shifted but keep $N_2$ fixed.  The associated partition will still have the same 
energy. But this cannot correspond to a minimal partition since the equal angle property does not hold; see also \cite{HHOT}. This argument extends to  arbitrary $m,n>0$ and $(\theta_1,\theta_2)$. 
\finbox  

\begin{remark}~\\
There exists  $0<b_0<1$ sufficiently close to $1$,  so that, for each irrational $b^2$ satisfying: $b_0<b< 1$,  only 
the first and the second eigenvalue together with their eigenfunctions are Courant sharp pairs.  
 This follows by counting. Remember $b<1$. The eigenvalues all have multiplicity $2$ 
or $4$. Suppose $n=0$ then $u_{m,0}$ has $ 2m$ nodal domains. So Courant sharpness can occur only for $\la_{m,0}=\la_{2m}$. 
This will not be the case if $|1-b|$ is small since then $\la_{0,n}$ will be eventually be below $\la_{m,0}$ hence 
$\la_{m,0}>\la_{2m}$. The case $m,n\ge 1$ has been treated above. 
\end{remark}

{\bf  Acknowledgements} \\
Thanks to P. de Soyres for his help for the pictures. The second author  had helpful discussions with Frank Morgan during the Dido conference in Carthage 2010.

{\scshape 
B. Helffer: D\'epartement de Math\'ematiques, Bat. 425,
Universit\'e Paris-Sud, 91 405 Orsay Cedex, France.

email: Bernard.Helffer@math.u-psud.fr\\

T. Hoffmann-Ostenhof: Department of Theoretical Chemistry, 1090 Wien,
W\"ahringerstrasse  17, Austria

email: thoffmann@tbi.univie.ac.at
}
\appendix
   %% \begin{figure}[h!]
  %%  \begin{center}
 %%  \includegraphics[height=20cm]{Documentgc.pdf}
  %%  \caption{Other $k$-partitions} \label{figder} 
  %%  \end{center}
  %%  \end{figure} Consider the torus $T(1,b)$ where $b>0$ and in addition $b^2$ is irrational. 
\section{Pictures}
\begin{figure}[h!]
\begin{center}
\caption{A $3$-partition of the torus without critical point.}\label{fig1} 
  \includegraphics[height=7cm]{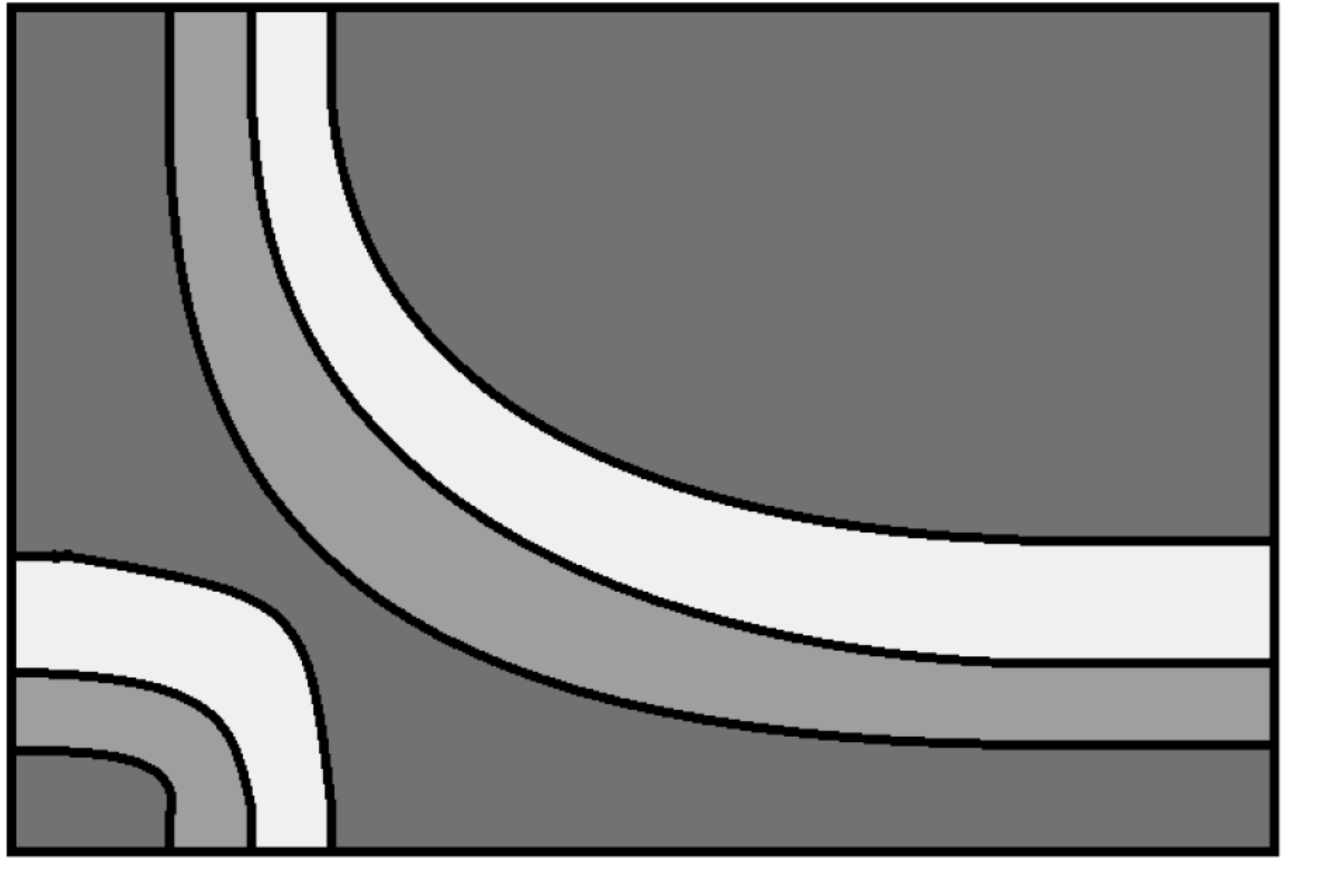}
   \caption{The lifted $3$-partition on the four-fold covering of the torus.}\label{fig2} 
 \includegraphics[height=11cm]{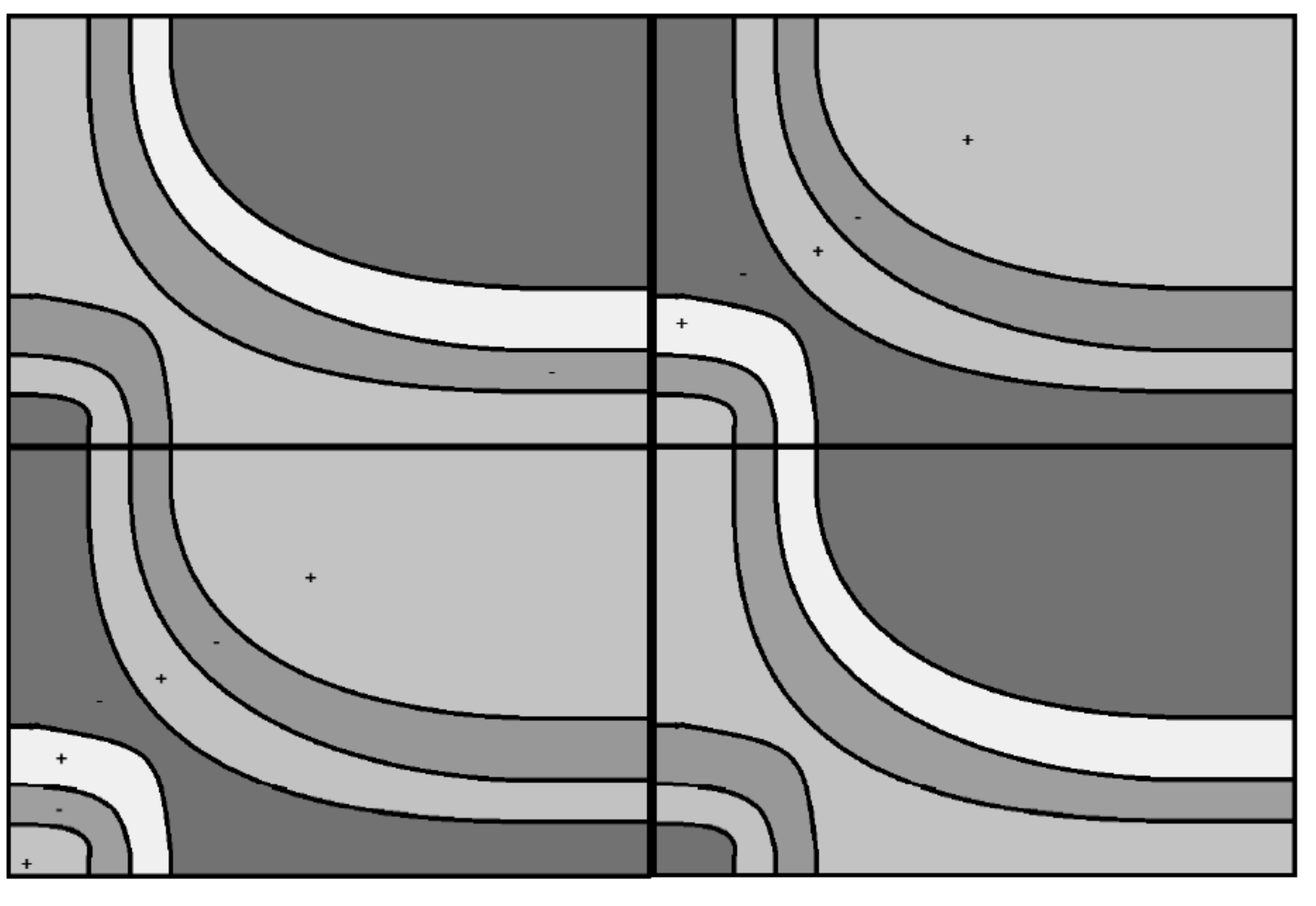}  \end{center}
\end{figure}
 \begin{figure}[h!]
   \begin{center}
  \caption{(p=3, q=2) and (p=4, q=2)}\label{figxa} 
  \includegraphics[height=15cm]{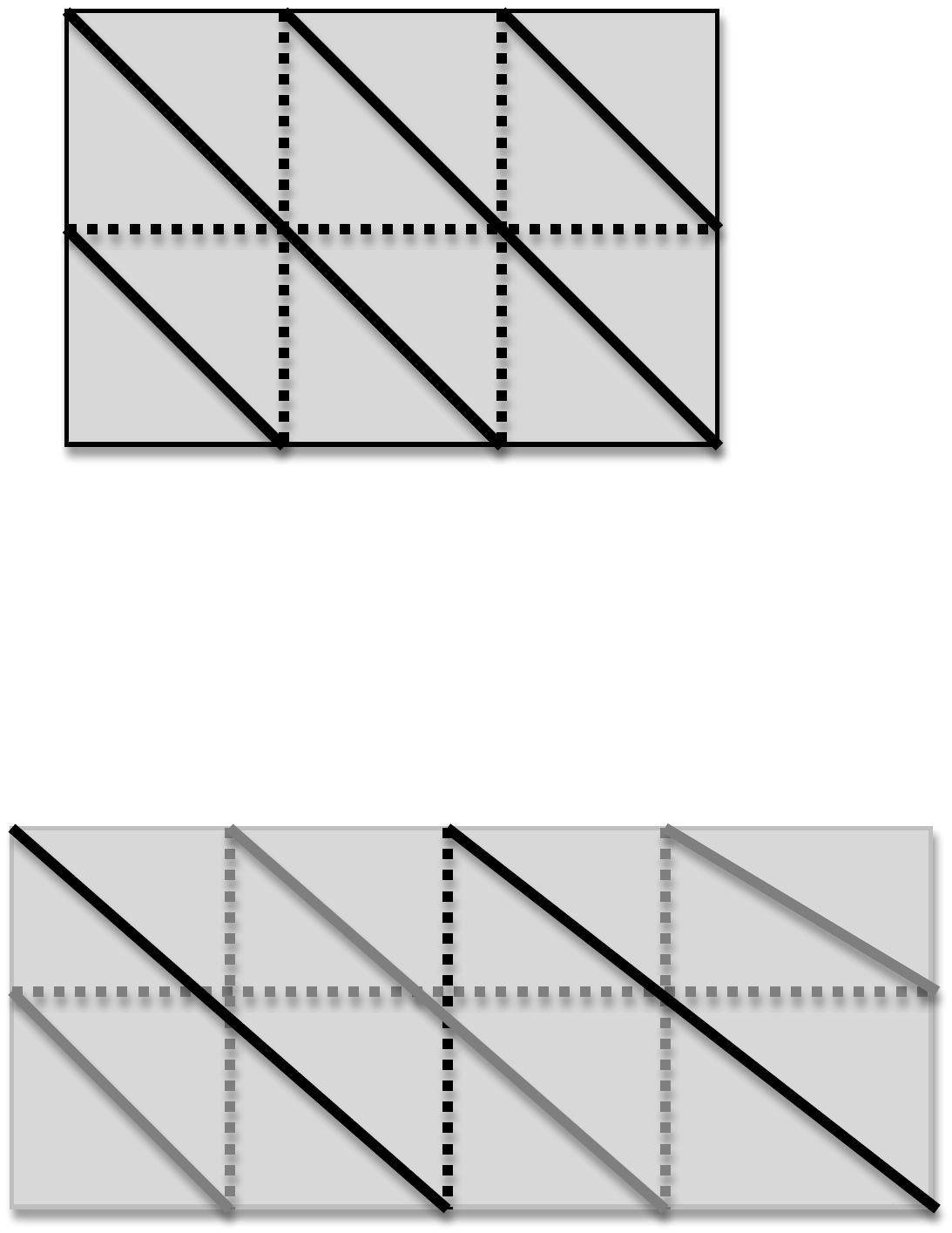}
  
 \end{center}
   \end{figure}

   \end{document}